\documentclass[a4paper]{scrartcl}
\usepackage{amsmath,amssymb,url,graphicx, verbatim}
\usepackage{varioref}
\usepackage{relsize}
\newtheorem{theorem}{Theorem}[section]
\newtheorem{lemma}[theorem]{Lemma}

\newtheorem{corollary}[theorem]{Corollary}

\DeclareMathOperator{\abs}{abs}
\numberwithin{equation}{section}
\newcommand{\extcite}[2]{\cite[#2]{#1}}

\begin{document}

\begin{center}
 \LARGE \textbf{PROBABILISTIC STAR DISCREPANCY BOUNDS FOR LACUNARY POINT SETS}
\end{center}

\vspace{3ex}

\begin{center}
 THOMAS L\"OBBE
\footnote{The author is supported by IRTG 1132}
\end{center}

\vspace{3ex}

\begin{abstract}
ABSTRACT. By a result of Heinrich, Novak, Wasilkowski and Wo\'zniakowski the inverse of the star discrepancy $n(d,\varepsilon)$ satisfies $n(d,\varepsilon)\leq c_{\abs}d\varepsilon^{-2}$. Equivalently for any $N$ and $d$ there exists a set of $N$ points in $[0,1)^d$ with star discrepacny bounded by $\sqrt{c_{\abs}\cdot d/N}$. They actually proved that a set of independent uniformly distributed random points satisfies this upper bound with positive probability. Although Aistleitner and Hofer later refined this result by proving a precise value of $c_{\abs}$ depending on the probability with which the inequality holds, so far there is no general construction for such a set of points known. In this paper we consider the sequence $(x_n)_{n\geq 1}=(\langle 2^{n-1}x_1\rangle)_{n\geq 1}$ for  a uniformly distributed point $x_1\in [0,1)^d$ and prove that the star discrepancy is bounded by $C\sqrt{d\log_2d/N}$. The precise value of $C$ depends on the probability with which this upper bound holds.
\end{abstract}

\section{Introduction}
A sequence of vectors $(x_n)_{n\geq 1}=(x_{n,1},\ldots,x_{n,d})_{n\geq 1}$ of real numbers in $[0,1)^d$ is called \textit{uniformly distributed modulo one} if
\begin{equation}
 \label{81}
 \lim_{N\to\infty}\frac{1}{N}\sum_{n=1}^N\mathbf{1}_{\mathcal{A}}(x_n)=\lambda(\mathcal{A})
\end{equation}
for any axis-parallel box $\mathcal{A}\subset [0,1)^d$ where $\mathbf{1}_{\mathcal{A}}$ denotes the indicator function on the set $\mathcal{A}$ and $\lambda$ denotes the Lebesgue-measure on $[0,1)^d$.
The 
star discrepancy of the first $N$ elements of $(x_n)_{n\geq 1}$ is defined by
\begin{equation}
 \label{65}
 \begin{aligned}
 D^*_N(x_1,\ldots,x_N)\:\:\: & =\:\:\: \sup_{\mathcal{A}\in\mathcal{B}^*}\left|\frac{1}{N}\sum_{n=1}^N\mathbf{1}_{\mathcal{A}}(x_n)-\lambda(\mathcal{A})\right|,
 \end{aligned}
\end{equation}
where 
$\mathcal{B}^*$ denotes the set of all axis-parallel boxes $\mathcal{A}=\prod_{i=1}^d[0,\beta_{i})\subset [0,1)^d$ with one corner in $0$. 
A sequence of points $(x_n)_{n\geq 1}$ is called a \textit{low-discrepancy sequence} if 
\begin{equation}
 \label{13}
 D_N^*(x_1,\ldots,x_N)\leq C\frac{(\log N)^d}{N}
\end{equation}
for all $N\geq 1$ and some absolute constant $C>0$. Furthermore Roth \cite{R80} showed that there exists a constant $C_d$ depending only on $d$ such that for any sequences $(x_n)_{n\geq 1}$ we have 
\begin{equation}
 \label{14}
 D_N^*(x_1,\ldots, x_N)\geq C_d\frac{(\log N)^{d/2}}{N}.
\end{equation}
Thus the asymptotic behaviour of a low-discrepancy sequence is not far from optimal. Nevertheless, if $N$ is small compared to $d$ then the upper bound in (\ref{13}) is not suitable. Therefore the \textit{inverse of the star discrepancy} was introduced. Let $n(d,\varepsilon)$ denote the smallest number $N$ such that there exists a $N$-element set of points in $[0,1)^d$ such that the star discrepancy of this point set is bounded by $\varepsilon$. In 2001 Heinrich, Novak, Wasilkowski and Wo\'zniakowski \cite{HNWW01} proved that
\begin{equation}
 \label{12}
 n(d,\varepsilon)\leq c_{\abs}d\varepsilon^{-2}
\end{equation}
holds for all $d\geq 1$ and $\varepsilon>0$ with some absolute constant $c_{\abs}>0$.  On the other hand Hinrichs \cite{H04} showed
\begin{equation}
 \label{10}
 n(d,\varepsilon)\geq c_{\abs}d\varepsilon^{-1}
\end{equation}
for all $d\geq 1$ and $\varepsilon>0$ and some possibly different absolute constant $c_{\abs}>0$. Thus the inverse of the star discrepancy depends linearly on the dimension, only the precise dependence on $\varepsilon$ is still unkown. By (\ref{9}) there exists a set of $N$ points in $[0,1)^d$ with
\begin{equation}
 \label{11}
 D^*_{N}(x_1,\ldots, x_N)\leq \sqrt{c_{\abs}}\sqrt{\frac{d}{N}}.
\end{equation}
In fact, Heinrich et \textit{al}. proved that a set of independent uniformly distributed random points, i.e. a Monte Carlo point set, satisfies (\ref{11}) with positive probability.  This result was later refined by Aistleitner and Hofer \cite{AH14} who gave an upper bound on $c_{\abs}$ depending on the probability with which (\ref{11}) is satisfied. Although they showed that even for moderate constants the inequality holds with high probability so far there is no general construction of a suitable point set known.

For a uniformly distributed point $x_1\in [0,1)^d$ let $(x_n)_{n\geq 1}$ be a sequence with $x_{n+1}=(x_{n+1,i})_{i=1,\ldots,d}=(\langle 2x_{n,i}\rangle)_{i=1,\ldots,d}$ for all $n\geq 1$ where $\langle \cdot\rangle$ denotes the fractional part of a rational number. Conze, Le Borgne and Roger \cite{CBR12} proved that a system of random variables $(f(x_n))_{n\geq 1}$ where $f:[0,1)^d\to\mathbb{R}$ is a centered indicator function on a box satisfies the Central Limit Theorem.  Thus the asymptotic behaviour of this sequence which is a particular example of a lacunary system $(f(M_nx)_{n\geq 1})$ which in general is defined by a centered one-periodic function $f$ with "nice" analytic properties and a fast growning sequence of $d\times d$ integer valued matrices satisfying a Hadamard gap condition
\begin{equation}
 \label{12a}
 || M_{n+k}^Tj||_{\infty}\geq q^k||M_n^T||_{\infty}
\end{equation}
for all $n,k\geq 1$, $j\in\mathbb{Z}^d$ with $0<\log_q||j||_{\infty}\leq  k$ and some absolute constant $q>1$
is similar to the behaviour of independent random variables.

The number of digits which are necessary to simulate $N$ points of this sequence with $H$ digits precision is of order $\mathcal{O}(d(H+N))$ and thus is much smaller than the number of digits to simulate $N$ independent random points which is $\mathcal{O}(dHN)$. Therefore we consider this randomized sequence $(x_{n})_{n\geq 1}$. We prove an upper bound on the star discrepancy which holds with high probability. Compared to (\ref{11}) this upper bound has up to some constant only an additional $\sqrt{\log_2d}$-factor. Our main result is stated in the following

\begin{theorem}
 \label{8}
 Let $N\geq 1$ and $d\geq 2$ be integers. 
 
 Then for any $0<\varepsilon<1$ the star discrepancy of the point set $(x_1,\ldots,x_N)$ satisfies
 \begin{equation*}
  D_N^*(x_1,\ldots,x_N)\leq (87-7d^{-1}\log\varepsilon)\sqrt{\frac{d\log_2d}{N}}
 \end{equation*}
 with probability at least $1-\varepsilon$.
\end{theorem}

\section{Preliminaries}

\begin{lemma}[Maximal Bernstein inequality, \extcite{EM96}{Lemma 2.2}]
 \label{191}
 For an integer $N\geq 1$ let $Z_1,\ldots,Z_N$ be a sequence of i.i.d. random variables with mean zero and variance $\sigma^2>0$ such that $|Z_1|\leq 1$. Then for any $t>0$ we have
 \begin{equation}
  \label{195}
   \mathbb{P}\left(\max_{M\in\{1,\ldots,N\}}\left|\sum_{n=1}^MZ_n\right|>t\right)\leq 2\exp\left(-\frac{t^2}{2N\sigma^2+2t/3}\right).
 \end{equation}
\end{lemma}

Let $v,w\in [0,1)^d$. We write $v\leq w$ if $v_i\leq w_i$ for all $i\in\{1,\ldots,d\}$. For some $\delta>0$ a set $\Delta$ of elements in $[0,1)^d\times [0,1)^d$ is called a $\delta$-bracketing cover if for every $x\in [0,1)^d$ there exists $(v,w)\in\Delta$ with $v\leq x\leq w$ and $\lambda(\overline{[v,w)})\leq\delta$ for $\overline{[v,w)}=[0,w)\backslash [0,v)$.
The following Lemma gives an upper bound on the cardinality of a $\delta$-bracketing cover.

\begin{lemma}[\extcite{G08}{Theorem 1.15}]
 \label{176}
 For any $d\geq 1$ and $\delta>0$ there exists some $\delta$-bracketing cover $\Delta$ with
 \begin{equation*}
  |\Delta|\leq \frac{1}{2}(2e)^d(\delta^{-1}+1)^d.
 \end{equation*}
\end{lemma}

\begin{corollary}
 \label{177}
 For any integers $d\geq 1$ and $h\geq1$ there exists a $2^{-h}$-bracketing cover $\Delta$ with
 \begin{equation*}
  |\Delta|\leq \frac{1}{2}(2e)^d(2^{h+2}+1)^d
 \end{equation*}
 such that for any $(v,w)\in\Delta$ and any $i\in \{1,\ldots,d\}$ we have
 \begin{eqnarray*}
  v_i & = & 2^{-(h+1+\lceil\log_2d\rceil)}a_i,\\
  w_i & = & 2^{-(h+2+\lceil\log_2d\rceil)}b_i
 \end{eqnarray*}
 for some integers $a_i\in\{0,1,\ldots,2^{h+1+\lceil\log_2d\rceil}\}$ and $b_i\in\{0,1,\ldots,2^{h+2+\lceil\log_2d\rceil}\}$.
\end{corollary}

\textit{Proof.} Let $\Delta$ be some $2^{-(h+2)}$-bracketing cover of $[0,1)^d$. By Lemma \ref{176} we have
\begin{equation*}
 |\Delta|\leq \frac{1}{2}(2e)^d(2^{(h+2)}+1)^d.
\end{equation*}
For $(v,w)\in\Delta$ and $i\in\{1,\ldots,d\}$ define
\begin{eqnarray*}
 y_{v,i} & = & \max\left\{2^{-(h+1+\lceil\log_2d\rceil)}a_i\leq v_i:a_i\in\mathbb{Z}\right\},\\
 z_{w,i} & = & \min\left\{2^{-(h+2+\lceil\log_2d\rceil)}b_i\geq w_i:b_i\in\mathbb{Z}\right\}.
\end{eqnarray*}
For $y_v=(y_{v,i})_{i\in\{1,\ldots,d\}}\in [0,1)^d$ we obtain
\begin{equation*}
 \lambda(\overline{[y_v,v)})\leq \sum_{i=1}^d2^{-(h+1+\lceil\log_2d\rceil)}\leq 2^{-(h+1)}.
\end{equation*}
Analogously for $z_w=(z_{w,i})_{i\in\{1,\ldots,d\}}\in [0,1)^d$ we have
\begin{equation*}
 \lambda(\overline{[z,z_w)})\leq 2^{-(h+2)}.
\end{equation*}
Thus we get
\begin{equation*}
 \lambda(\overline{[y_v,z_w)})\leq \lambda(\overline{[y_v,v)})+\lambda(\overline{[v,w)})+\lambda(\overline{[w,z_w)})\leq 2^{-h}.
\end{equation*}
Set $\tilde{\Delta}=\{(y_v,z_w):(v,w)\in\Delta\}$. Since $\Delta$ is a $2^{-(h+2)}$-bracketing cover for any $x\in[0,1)^d$ there exists $(v,w)\in\Delta$ and $(y_v,z_w)\in\tilde{\Delta}$ with $y_v\leq v\leq x\leq w\leq z_w$. Therefore $\tilde{\Delta}$ is a $2^{-h}$-bracketing cover and the conclusion of the proof follows by $|\tilde{\Delta}|\leq |\Delta|$.

\section{Proof of main theorem}

The proof of this Theorem is mainly based on \cite{A13}. For some integers $N\geq 1$ and $d\geq 1$ we simply write
\begin{equation*}
 D^d_N(x_{n,i})=D^d_N((x_{1,1},\ldots,x_{1,d}),\ldots,(x_{N,1},\ldots,x_{N,d})).
\end{equation*}

For $N\geq 1 $ and $d\geq 1$ set
\begin{equation}
 \label{180}
 H=\left\lceil\frac{\log_2N}{2}-\frac{\log_2(d\log_2d)}{2}-2\right\rceil.
\end{equation}
As a consequence for any $h\in\{0,\ldots,H\}$ we have
\begin{equation}
 \label{181}
 \sqrt{d\log_2d}\sqrt{N}\leq 2^{-h}N.
\end{equation}
For any $h\in\{1,\ldots,H\}$ let $\Delta_h$ be a $2^{-h}$-bracketing cover of $[0,1)^d$. By Corollary \ref{177} we may assume
\begin{equation}
 |\Delta_h|\leq \frac{1}{2}(2e)^d(2^{h+2}+1)^d.
\end{equation}
For any $y\in [0,1)^d$ we now define a finite sequence of points $\beta_h(y)$ for $h\in\{0,\ldots,H+1\}$ in the following manner. Let $(v,w)\in\Delta_H$ be such that $v\leq y\leq w$. We set $\beta_{H+1}(y)=w$ and $\beta_H(y)=v$.
The points $\beta_1(y),\ldots,\beta_{H-1}(y)$ are defined by induction. Thus assume that for some $h\in\{1,\ldots,H-1\}$ the point $\beta_{h+1}(y)$ is already defined. Let $(v,w)\in\Delta_{h}$ with $v\leq \beta_{h+1}(y)\leq w$ and set $\beta_{h}(y)=v$.
Moreover set $\beta_0(y)=0$. Therefore we observe
\begin{equation*}
 0=\beta_0(y)\leq \beta_1(y)\leq \cdots \leq \beta_H(y)\leq x\leq \beta_{H+1}(y)\leq 1.
\end{equation*}
For $h\in\{0,\ldots,H-1\}$ we have $(\beta_h(y),w)\in\Delta_h$ for some point $w\in [0,1)^d$.
Furthermore we have $(\beta_H(y),\beta_{H+1}(y))\in\Delta_H$.
Then by Corollary \ref{177} for $h\in\{0,\ldots,H+1\}$ and $i\in\{1,\ldots, d\}$ there exist integers $a_{h,i}\in\{0,\ldots,2^{h+1+\log_2d}\}$ such that
\begin{equation}
 \label{182}
 (\beta_h(y))_i=2^{-(h+1+\log_2d)}a_{h,i}.
\end{equation}
For $h\in\{0,\ldots,H\}$ set $K_h(y)=\overline{[\beta_h(y),\beta_{h+1}(y))}$. Note that the sets $K_h(y)$ are pairwise disjoint and satisfy
\begin{equation}
 \label{183}
 \bigcup_{h=0}^{H-1}K_h(x)\subseteq [0,x)\subseteq\bigcup_{h=0}^{H}K_h(x)
\end{equation}
By definition $\beta_h(y)\leq \beta_{h+1}(y)\leq w$ for some $w\in [0,1)^d$ with $(\beta_h(y),w)\in\Delta_h$ and hence
\begin{equation}
 \label{184}
 \lambda(K_h(y))\leq \lambda\left(\overline{[\beta_h(y),w)}\right)\leq 2^{-h}
\end{equation}
for any $h\in\{0,\ldots,H\}$.
Now define 
\begin{equation*}
S_h=\left\{\overline{[\beta_h(y),\beta_{h+1}(y))}:y\in [0,1)^d\right\}. 
\end{equation*}
Observe that we may define the points $\beta_h$ such that $\beta_h(y)=\beta_h(z)$ for $y,z\in [0,1)^d$ with $\beta_{h+1}(y)=\beta_{h+1}(z)$. Therefore by Corollary \ref{177} we have
\begin{equation}
 \label{185}
 |S_h|=\left|\left\{\beta_{h+1}(y):y\in [0,1)^d\right\}\right|\leq |\Delta_{h+1}|\leq \frac{1}{2}(2e)^d(\sqrt{5})^{(h+3)d}
\end{equation}
for any integer $h\in\{0,\ldots,H\}$.
Note that hereafter we skip the point $y$ in the notation of the points $\beta_h$ and the sets $K_h$ to simplify notations.
Then by (\ref{183}) we have
\begin{equation}
 \label{186}
  \sum_{n=1}^N\mathbf{1}_{[0,y)}(x_n) \geq  \sum_{n=1}^{N}\mathbf{1}_{[0,\beta_H)}(x_n) = \sum_{h=0}^{H-1}\sum_{n=1}^N\left(\mathbf{1}_{K_{h}}(x_n)-\lambda(K_h)\right).
\end{equation}
Analogously we also get
\begin{equation}
 \label{187}
  \sum_{n=1}^N\mathbf{1}_{[0,y)}(x_n) \leq  \sum_{n=1}^{N}\mathbf{1}_{[0,\beta_{H+1})}(x_n) = \sum_{h=0}^{H}\sum_{n=1}^N\left(\mathbf{1}_{K_{h}}(x_n)-\lambda(K_h)\right).
\end{equation}
By using Bernstein inequality we now shall give a lower bound on the probability that the inequality
\begin{equation}
 \label{204}
  \left|\sum_{n=1}^N\left(\mathbf{1}_{K_h}(x_n)-\lambda(K_h)\right)\right|> t
\end{equation}
holds simultaneously for all $h\in\{0,\ldots,H\}$ and some $t>0$ to specified later . 
Observe that in general the random variables $f_{K_h}(x_n)=\mathbf{1}_{K_h}(x_n)-\lambda(K_h)$ are not independent. Thus we may not apply the Bernstein inequality directly. Therefore we  decompose the set of numbers $\{1,\ldots,N\}$ into several modulo classes. If the distance between two consecutive indices $n_l,n_{l+1}$ in the same class is large enough, i.e.  $n_{l+1}-n_l\geq h+2+\lceil\log_2d\rceil$, the random variables are stochastically independent, i.e.
\begin{equation}
 \label{194}
 \mathbb{P}\left(f_{K_h}(x_{n_1})=c_1,\ldots,f_{K_h}(x_{n_k})=c_k\right)=\prod_{l=1}^k\mathbb{P}\left(f_{K_h}(x_{n_l})=c_l\right).
\end{equation}
We only prove the case $k=2$. The general case follows by induction.

By (\ref{182}) the set $K_h$ is a union of axis-parallel boxes such that each corner of any box is of the form
\begin{equation}
 \label{193}
 \left(2^{-(h+2+\lceil\log_2d\rceil)}a_{1},\ldots,2^{-(h+2+\lceil\log_2d\rceil)}a_{d}\right)
\end{equation}
such that $a_i\in\{0,1,\ldots,2^{h+2+\lceil\log_2d\rceil}\}$ for any $i\in\{s+1,\ldots,d\}$.
Furthermore let $n,n'\in\{1,\ldots,N\}$ be two indices with $n'-n\geq h+2+\lceil\log_2d\rceil$. We define a decomposition of $[0,1)^{d}$ by
\begin{multline*}
 \Sigma=\left\{\prod_{i=1}^d\left[2^{-(n'-1+\lceil\log_2d\rceil)}a_i,2^{-(n'-1+\lceil\log_2d\rceil)}(a_i+1)\right):\right.\\
 \left.a_i\in\left\{0,1,\ldots,2^{n'-1+\lceil\log_2d\rceil}-1\right\},i\in\{1,\ldots,d\}\right\}.
\end{multline*}
Note that by (\ref{193}) the function $f_{K_h}$ is constant on any box $\mathcal{B}\in\Sigma$. For some $c_1\in\mathbb{R}$ define
\begin{equation*}
 \Sigma_{c_1}=\left\{\mathcal{B}\in\Sigma:f_{K_h}(x_{n})=c_1 \textnormal{ for all }x_1=(x_{1,1},\ldots,r_{1,d})\in\mathcal{B}\right\}.
\end{equation*}
Since $x_{n',i}=2^{n'-1}x_{1,i}$ for all $i\in\{1,\ldots,d\}$ we have $f_{K_h}(x_{n'})=f_{K_h}(x'_{n'})$ where $x'_{n'}=(x'_{n',1},\ldots,x'_{n',d})$ with $x'_{n,i}=2^{n'-1}x'_{1,i}$ is an
instance of the matrix for some initial value $x'_1=(x'_{1,1},\ldots,x'_{1,d})$ with $x'_{1,i}=x_{1,i}+2^{-(n'-1+\lceil\log_2d\rceil)}a_i$ and $a_i\in\{0,1,\ldots,2^{n'-1+\lceil\log_2d\rceil}-1\}$ for all $i\in\{1,\ldots,d\}$.
Therefore for any $c_2\in\mathbb{R}$ and any $\mathcal{B},\mathcal{B}'\in\Sigma$ we have
\begin{equation*}
 \mathbb{P}\left(f_{K_h}(x_{n'})=c_2|x_1\in\mathcal{B}\right)=\mathbb{P}\left(f_{K_h}(x_{n'})=c_2|x_1\in\mathcal{B}'\right).
\end{equation*}
Hence for any $c_2\in\mathbb{R}$ and any $\mathcal{B}\in\Sigma$ we get
\begin{eqnarray*}
 \mathbb{P}\left(f_{K_h}(x_{n'})=c_2\right) & = & \sum_{\mathcal{B}'\in\Sigma}\mathbb{P}\left(f_{K_h}(x_{n'})=c_2|r_1\in\mathcal{B}'\right)\mathbb{P}(x_1\in\mathcal{B}')\\
 & = & \mathbb{P}\left(f_{K_h}(x_{n'})=c_2|x_1\in\mathcal{B}\right)\sum_{\mathcal{B}'\in\Sigma}\mathbb{P}(x_1\in\mathcal{B}')\\
 & = & \mathbb{P}\left(f_{K_h}(x_{n'})=c_2|x_1\in\mathcal{B}\right).
\end{eqnarray*}
Moreover for any $c_1,c_2\in\mathbb{R}$ we obtain
\begin{multline*}
 \mathbb{P}\left(f_{K_h}(r_{n'})=c_2|f_{K_h}(x_n)=c_1\right)\\
 \begin{aligned}
  = & \frac{\mathbb{P}\left(f_{K_h}(x_{n'})=c_2,f_{K_h}(x_n)=c_1\right)}{\mathbb{P}\left(f_{K_h}(x_n)=c_1\right)}\\
  = & \frac{\sum_{\mathcal{B}\in\Sigma}\mathbb{P}\left(f_{K_h}(x_{n'})=c_2,f_{K_h}(x_n)=c_1|x_1\in\mathcal{B}\right)\mathbb{P}(x_1\in\mathcal{B})}{\mathbb{P}\left(f_{K_h}(x_n)=c_1\right)}\\
  = & \sum_{\mathcal{B}\in\Sigma_{c_1}}\mathbb{P}\left(f_{K_h}(x_{n'})=c_2|x_1\in\mathcal{B}\right)\frac{\mathbb{P}(x_1\in\mathcal{B})}{\mathbb{P}\left(f_{K_h}(x_n)=c_1\right)}\\
  = & \mathbb{P}\left(f_{K_h}(x_{n'})=c_2\right).
 \end{aligned}
\end{multline*}
Thus (\ref{194}) is proved.
Set $\kappa=\kappa_h=\lceil \log_2(h+2+\lceil\log_2d\rceil)\rceil$. Furthermore set $Q(N,\kappa,\gamma)=\{n\in\{1,\ldots,N\}:n\equiv\gamma (\mod 2^{\kappa})\}.$

Then for $h\in \{0,\ldots, H\}$ by Lemma \ref{191} we have
\begin{multline}
 \label{1}
 \mathbb{P}\left(\left|\sum_{n=1}^N\left(\mathbf{1}_{K_h}(x_n)-\lambda(K_h)\right)\right|>t\right)\\
 \begin{aligned}
  \leq & \sum_{\gamma=1}^{2^{\kappa}}\mathbb{P}\left(\left|\sum_{n\in Q(N,\kappa,\gamma)}\mathbf{1}_{K_h}(r_n)-\lambda(K_h)\right|>\frac{t}{2^{\kappa}}\right)\\
  \leq & 2\sum_{\gamma=1}^{2^{\kappa}}\exp\left(-\frac{t^2/2^{2\kappa}}{2\left(\sum_{n\in Q(N,\kappa,\gamma)}1\right)\lambda(K_h)(1-\lambda(K_h))+2t/(3\cdot 2^{\kappa})}\right)\\
  \leq & 2^{\kappa+1}\exp\left(\frac{t^2/2^{\kappa}}{4N\cdot 2^{-h}+2t/3}\right).
 \end{aligned}
\end{multline}
For $h\geq 1$ set $t=C_1\sqrt{d\log_2d}\sqrt{N}\sqrt{h\cdot 2^{-h}}$ for a constant $C_1>0$ to be specified later. By (\ref{181}) we observe
\begin{equation*}
 \frac{2t}{3}\leq \frac{2}{\sqrt{3}}C_1\log_2d\cdot N\cdot 2^{-h}.
\end{equation*}
Thus we get
\begin{equation*}
 \frac{t^2}{4\cdot 2^{-h}N+2t/3} \geq \frac{C_1^2d\log_2d\cdot 2^{-h}hN}{4\cdot 2^{-h}N+2/\sqrt{3}\cdot C_1\log_2d\cdot 2^{-h}N} \geq\frac{C_1^2dh}{4+2/\sqrt{3}\cdot C_1}.
\end{equation*}
Plugging this into (\ref{1}) we obtain
\begin{multline}
 \label{2}
 \mathbb{P}\left(\left|\sum_{n=1}^N\left(\mathbf{1}_{K_h}(x_n)-\lambda(K_h)\right)\right|>t\right)\\
 \begin{aligned}
  \leq & 2\exp\left(\kappa\log 2-\frac{C_1^2}{4+2/\sqrt{3}\cdot C_1}dh\right)\\
  \leq & 2\exp\left(\lceil \log_2(h+2+\lceil\log_2d\rceil)\rceil\log 2-\frac{C_1^2}{4+2/\sqrt{3}\cdot C_1}dh\right)\\
  \leq & 2\exp\left(-\left(\frac{C_1^2}{4+2/\sqrt{3}\cdot C_1}-1\right)dh\right)\\
 \end{aligned}
\end{multline}
where the last inequality holds for $d\geq 2,h\geq 1$. Similarly for $h=0$ we set $t=C_2\sqrt{d\log_2d}\sqrt{N}$ for a constant $C_2>0$ to be specified later. Repeating the above calculation we show
\begin{equation}
 \label{3}
 \mathbb{P}\left(\left|\sum_{n=1}^N\left(\mathbf{1}_{K_h}(x_n)-\lambda(K_h)\right)\right|>t\right)\leq 2\exp\left(-\left(\frac{C_2^2}{4+2/3\cdot C_2}-1\right)d\right)
\end{equation}
for $d\geq 2$.
Now define
\begin{equation}
 \label{4}
 C_3=\frac{C_1^2}{4+2/\sqrt{3}\cdot C_1}-1, \quad C_4=\frac{C_2^2}{4+2/3\cdot C_2}-1.
\end{equation}
By (\ref{185}) the statement of the Theorem immediately follows if we show
\begin{equation}
 \label{5}
 1-\frac{1}{2}(2e)^d(\sqrt{5})^{3d}\cdot 2e^{-C_4d}-\frac{1}{2}(2e)^d\sum_{h=1}^H(\sqrt{5})^{(h+3)d}\cdot 2e^{-C_3dh}\geq 1-\varepsilon.
\end{equation}
Thus it is enough to choose constants $C_3,C_4$ large enough such that
\begin{equation}
 \label{6}
 \frac{1}{2}(2e)^d(\sqrt{5})^{3d}\cdot 2e^{-C_4d}\leq \frac{\varepsilon}{2}
\end{equation}
and
\begin{equation}
 \label{7}
 \frac{1}{2}(2e)^d(\sqrt{5})^{(h+3)d}\cdot 2e^{-C_3dh}\leq \frac{\varepsilon}{2^{h+1}}
\end{equation}
for all $h\in\{1,\ldots,H\}$.
Observe that (\ref{6}) is satisfied for
\begin{equation}
 \label{8a}
 C_4=4.46-\frac{\log\varepsilon}{d} \geq 1+\log 2+1.5\log 5 +\frac{\log 2}{d}-\frac{\log\varepsilon}{d}.
\end{equation}
Similarly (\ref{7}) is equivalent to
\begin{equation*}
 (1+2\log 2+1.5\log 5)d+\frac{\log 5}{2}dh+\log2 \cdot h-\log\varepsilon \leq C_3dh.
\end{equation*}
Since $h\geq 1$ we may choose
\begin{equation}
 \label{9}
 C_3=6.31-d^{-1}h^{-1}\log\varepsilon.
\end{equation}
By (\ref{4}) we may set
\begin{eqnarray}
 C_1 & = & 15.465-1.155d^{-1}\log\varepsilon,\\
 C_2 & = & 9.864-2/3\cdot d^{-1}\log\varepsilon.
\end{eqnarray}
Thus with probability at least $1-\varepsilon$ by (\ref{186}), (\ref{2}) and (\ref{3}) we have
\begin{multline*}
 \sum_{n=1}^N\mathbf{1}_{[0,y)}(x_n)\\
 \begin{aligned}
  \leq & \sum_{n=1}^N\sum_{h=1}^{H+1}\mathbf{1}_{K_h}(x_n)\\
  \leq & \sum_{n=1}^N\left(\lambda([0,\beta_1))+\left(9.864-\frac{2\log\varepsilon}{3d}\right)\sqrt{\frac{d\log_2d}{N}}\right)\\
           &  +\sum_{n=1}^N\sum_{h=1}^H\left(\lambda(\overline{[\beta_h,\beta_{h+1})})+\left(15.465-\frac{1.155\log\varepsilon}{d}\right)\sqrt{\frac{d\log_2d}{N}}\sqrt{2^{-h}h}\right)\\
  \leq & \sum_{n=1}^N\left(\lambda([0,y))+\lambda([y,\beta_{H+1}))+\left(82.357-6.081d^{-1}\log\varepsilon\right)\sqrt{\frac{d\log_2d}{N}}\right)\\
  \leq & \sum_{n=1}^N\lambda([0,y))+\sum_{n=1}^N\left(86.357-6.081d^{-1}\log\varepsilon\right)\sqrt{\frac{d\log_2d}{N}}.
 \end{aligned}
\end{multline*}
Analogously we obtain
\begin{equation*}
 \sum_{n=1}^N\mathbf{1}_{[0,y)}(x_n)\geq\sum_{n=1}^N\lambda([0,y))-\sum_{n=1}^N\left(86.357-6.081d^{-1}\log\varepsilon\right)\sqrt{\frac{d\log_2d}{N}}
\end{equation*}
with probability at least $1-\varepsilon$.

\vspace{3ex}

\small{DEPT. OF MATHEMATICS, BIELEFELD UNIV., P.O.Box 100131, 33501 Bielefeld, Germany\\
\textit{E-Mail address:} \url{tloebbe@math.uni-bielefeld.de}}


\begin{thebibliography}{99999}
\bibitem{A13} Aistleitner, C.: On the inverse of the discrepancy for infinite dimensional infinite sequences, J. Complexity 29, 182-194 (2013)
\bibitem{AH14} Aistleitner, C., Hofer, M.: Probabilistic discrepancy bound for Monte Carlo point sets, Math. Comp. 83., 1373-1381 (2014)
\bibitem{CBR12} Conze, J.-P., Le Borgne, S., Roger, M.: Central limit theorem for stationary products of toral automorphisms, Discrete Contin. Dyn. Syst. 32, 1597-1626 (2012)
\bibitem{EM96} Einmahl, U., Mason, D.M.: Some universal results on the behavior of increments of partial sums, Ann. Prob. 24, 1388-1407 (1996)
\bibitem{G08} Gnewuch, M.: Bracketing numbers for axis-parallel boxes and applications to geometric discrepancy, J. Complexity 24, 154-172 (2008)
\bibitem{HNWW01} Heinrich, S., Novak E., Wasilkowski, G.W., Wo\'zniakowski, H.: The inverse of the star-discrepancy depends linearly on the dimension, Acta Arith. 96, 279-302 (2001)
\bibitem{H04} Hinrichs, A.: Covering numbers, Vapnik-\v Cervonenkis classes and bounds for the star-discrepancy, J. Complexity 20, 477-483 (2004)
\bibitem{R80} Roth, K.F.: On irregularities of distribution I-IV, Mathematika 1, 73-79 (1954), Comm. Pure Appl. Math. 29, 739-744 (1976), Acta Arith. 35, 373-384 (1979) and Acta Arith. 37, 67-75 (1980)
\end{thebibliography}
\end{document}